\documentclass[10pt]{article}
\usepackage{amsfonts}
\usepackage{amssymb}
\usepackage{latexsym}
\usepackage{mathptm, times}

\title{On Lattice Barycentric Tetrahedra}
\author{Brian Mazur}

\newcommand{\Z}{\mathbb{Z}}

\renewcommand{\labelitemi}{$\Diamond$}

        {
        \smallskip      
        \stepcounter{theorem}
        \noindent
        {\bf Example \arabic{section}.\arabic{theorem}:} \ \ }
        {\hspace*{\fill}{$\Diamond$}
        \smallskip}

\newenvironment{proof}[1][]
        {
        \noindent
        {\bf Proof{#1}:  }
        }
        {\hspace*{\fill}{$\Box$}\smallskip}

        {
        \noindent
        {\bf Sketch{#1}:  }
        }
        {\hspace*{\fill}{$\Box$}\smallskip}

        {
        \noindent
        {\bf Reference:  }
        }
        {\hspace*{\fill}{$\odot$}\smallskip}

\newtheorem{theorem}{Theorem}[section]
\newtheorem{proposition}[theorem]{Proposition}
\newtheorem{lemma}[theorem]{Lemma}
\newtheorem{definition}[theorem]{Definition}

\begin{document}

\maketitle
\begin{abstract}
\noindent
A lattice tetrahedron $T \subset \mathbb{R}^3$
is a tetrahedron whose four vertices are
all in the lattice $\mathbb{Z}^3$.  Lattice tetrahedra are preserved by those affine linear maps 
of the form \(\vec{v} \mapsto A\vec{v} + \vec{b}\), such that \(A\) is an element of \(GL(3,{\Z}) \) and 
\begin{math} \vec{b} \end{math} is an element of the lattice \begin{math} \mathbb{Z}^3 \end{math}.  Such affine linear
maps are called unimodular maps.  We say that a lattice tetrahedron whose barycentre is its only non-vertex lattice  point is
lattice barycentric. 
The notation \begin{math}T(a,b,c)\end{math} describes that lattice tetrahedron with vertices 
\(
\{{0},{e}_1,{e}_2, a{e}_1+b{e}_2+
c{e}_3\}.
\)
Our result is then that all such lattice barycentric tetrahedra  \begin{math}T(a,b,c)\end{math} are unimodularly equivalent
to \begin{math}T(3,3,4)\end{math} or \begin{math} T(7,11,20)\end{math}.
\end{abstract}

\tableofcontents

\newpage

\section{Introduction}

\begin{definition}
\label{def:fp}
For $v_1, v_2, v_3, v_4 \in \mathbb{R}^3$, the tetrahedron 
$ 
T = T(v_1, v_2, v_3 ,v_4) 
$
is
\begin{equation}
T(v_1, v_2, v_3 ,v_4) =
\bigg\{ v = \sum_{i=1}^4 \lambda_i v_i \ :\   0\leq\lambda_i
\leq 1, \sum_{i=1}^4 \lambda_i = 1\bigg \}
\end{equation}

\smallskip
\noindent
A lattice tetrahedron $T$  is a tetrahedron such that
\(
\{v_1, v_2, v_3 , v_4\} \subset {\Z}^3,\; \mbox{where}\; {\Z}^3 = \{ a{e_1}+b{e_2}+c{e_3} \ : \ a,b,c \in {\Z} \}
\).   

\smallskip
\noindent 
A lattice tetrahedron 
\(T\) is said
to be fundamental if \(\partial T \; \cap \; {\Z}^3 = 
\{v_1, v_2 , v_3, v_4\}\). 

\smallskip
\noindent
By \(int(T)\) we mean those points \(v\) with \( 0 < \lambda_i < 1, 1 \leq i \leq 4 \).
A lattice tetrahedron \begin{math} T \end{math} is said to be primitive if it is fundamental and \begin{math}
\mbox{int} (T) \; \cap \; {\Z}^3 = \oslash \end{math}.

\end{definition}

\noindent
The barycentre, or centroid, of any tetrahedron \begin{math} T(v_1, v_2, v_3 ,v_4) \end{math} is the point \begin{math} BC[T(v_1, v_2, v_3 ,v_4)]
\end{math}
formed by an equal weighting of the vertices, that is,

\begin{equation}BC[T(v_1, v_2, v_3 ,v_4)] = \frac{1} {4} v_1 + \frac{1} {4} v_2 + \frac {1} {4} v_3 + \frac {1} {4} v_4 \end{equation}

\begin{definition}
A fundamental \begin{math} T(v_1, v_2, v_3, v_4)\end{math} with \begin{math} \{v_1, v_2, v_3, v_4\} \subset \mathbb{Z}^3 \end{math} and \(
int(T) \; \cap \;  \mathbb{Z}^3 = \{BC[T(v_1, v_2, v_3 ,v_4)]\} \)
is called lattice barycentric.

\end{definition} 

For unimodular equivalance we will work with elements of \( GL(3, \mathbb{Z}) \), which is defind as follows:

\begin{equation} GL(3,\mathbb{Z}) = \bigg\{ A = (A_{ij})\; \mbox{a} \; 3 \times 3 \; \mbox{matrix}   \ : \ A_{ij} \in \mathbb{Z} \ ; \  det(A) = \pm 1 
\bigg\}.
\end{equation}

\noindent
We say two tetrahedra \begin{math} T_1, T_2\end{math}
are unimodularly equivalent if there exists an
\begin{math} A \in GL(3,{\Z})\end{math} and
\begin{math} \vec{b} \in {\Z}^3\end{math} such that the map 
\begin{math} \vec{v}
\mapsto  A\vec{v}+\vec{b}\end{math} carries \begin{math}T_1\end{math} 
bijectively to 
\begin{math}T_2\end{math}.  Such 
\(
\vec{v} \mapsto A \vec{v} + \vec{b}
\)
are called (affine) unimodular maps.
These maps are precisely the maps which are lattice and volume preserving.  All lattice tetrahedra are unimodularly equivalent to
some $T(a,b,c) = T(0, e_1, e_2, a{e_1} + b{e_2} + c{e_3})$\cite[Thm. 5.2]{reznick:86a}.  We now state the main result.  

\begin{proposition}\label{Prop:main}
Every lattice barycentric tetrahedron $T$ is unimodularly equivalent to either $T(3,3,4)$ or $T(7,11,20)$.
\end{proposition}

\section{Outline Of The Proof And Basic Notation}
\subsection{Tetrahedra}

\begin{definition}
A \emph{grounded tetrahedron} $T$ has vertices $0,  e_1, e_2,$ and $ x e_1 + y e_2 + z e_3$. 
Such grounded tetrahedra are denoted $T(x,y,z)$.
\end{definition} 

We start the proof by placing certain greatest common divisor conditions (GCD) on 
$(a,b,c)$  such that $T(a,b,c)$ is fundamental per Definition \ref{def:fp}.  Second we recall further GCD
conditions requiring that $T(a,b,c)$ be primitive.

\begin{itemize}
\item[\labelitemi]\label{Fact:f}Fact 1: {$T(a,b,c)$ is fundamental if and only if $GCD(a,c) = GCD(b,c) = GCD(a+b-1, c) = 1$.}\cite[Thm. 5.2]{reznick:86a}
\item[\labelitemi]\label{Fact:p}Fact 2: {$T(a,b,c)$ is primitive if and only if $T(a,b,c)$ is fundamental, and also $a \equiv 1 \bmod c$ or 
$b \equiv 1 \bmod c$ or $a + b \equiv 0 \bmod c$.} \cite[Thm. 5.5 (Reeve-White-Howe-Scarf)]{reznick:86a} 
\end{itemize}

Introducing new notation, put the barycentre, $BC[T(a,b,c)] = (\alpha, \beta, \gamma)$, where $\alpha = \frac{a+1}{4}$, $\beta = \frac{b+1}{4}$, 
and $\gamma = \frac{c}{4}$.  For lattice barycentric tetrahedron, $(\alpha, \beta, \gamma)$ is required to be a lattice point.  This forces 
$a \equiv 3 \bmod 4$, $ b \equiv 3 \bmod 4$,  and $c \equiv 0 \bmod 4$.
Using the barycentre as a common vertex, one may cone over the triangular faces of the tetrahedra $T$ to produce four
sub-tetrahedra, $\mathcal{T}_j, j = 1,2,3,4$.  We would then like to analyse these sub-tetrahedra for primitiveness.  For if we know each sub-tetrahedra
$\mathcal{T}_j$ is primitive, then we know that $T(a,b,c)$ is lattice barycentric.  The conditions for primitiveness, however, may only be 
applied on the grounded sub-tetrahedron, that is, the unique sub-tetrahedron with vertices $0, e_1, e_2, BC[T(a,b,c)]$.  We define this to be 
$\mathcal{T}_4$ of $T$.  One may find unimodular maps which bring a sub-tetrahedra 
$\mathcal{T}_j$ of $T$ into the ground position.  These maps move the entire tetrahedron $T$ such that the apex, $(a,b,c)$, is sent to another point 
$(\tilde{a}, \tilde{b}, \tilde{c}) \in {\Z}^3$.  The primitivity conditions may then be applied to $T(\tilde{a}, \tilde{b}, \tilde{c})$ for 
$\tilde{\alpha} = \frac{\tilde{a} + 1}{4}, \tilde{\beta} = \frac{\tilde{b} + 1}{4},$ and $\tilde{\gamma} = \gamma = \frac{\tilde{c}}{4} = \frac{c}{4}$.

%\setlength{\unitlength}{1cm}
%
%\begin{figure}[t]
%\begin{center}
%\begin{picture}(6,6)
%
%\put(0,0){\line(1,1){2.0}}
%\put(2,2){\line(1,0){4.0}}
%\put(0,0){\line(3,1){6.0}}
%\put(0,0){\line(1,3){2.0}}
%\put(6,2){\line(-1,1){4.0}}
%\put(2,6){\line(0,-1){4.0}}
%\put(0,0){\line(3,4){3.0}}
%\put(2,2){\line(1,2){1.0}}
%\put(3,4){\line(-1,2){1.0}}
%\put(6,2){\line(-3,2){3.0}} 
%
%\end{picture}
%\end{center}
%\caption{This is an example of a tetrahedron which has been coned into its four sub-tetrahedra.}
%\end{figure}  

\subsection{Construction Of Unimodular Maps}
We construct unimodular maps, $\mathbf{h}_j(v)$, which carry the respective sub-tetrahedra into the grounded position as follows:

\begin{equation}
\label{Eq:maps}
\begin{array} {l}
 \mathbf h_1(v)  \\
0  \mapsto  (\tilde{a}, \tilde{b}, \tilde{c})\\
e_1  \mapsto  e_1\\
e_2  \mapsto  e_2\\
(a,b,c)  \mapsto  0\\
\end{array}
\begin{array} {l}
 \mathbf h_2(v)  \\
0  \mapsto  0\\
e_1  \mapsto  (\tilde{a}, \tilde{b}, \tilde{c}) \\
e_2  \mapsto  e_2\\
(a,b,c)  \mapsto  e_1\\
\end{array}
\begin{array} {l}
 \mathbf h_3(v)  \\
0  \mapsto  0\\
e_1  \mapsto  e_1\\
e_2  \mapsto  (\tilde{a}, \tilde{b}, \tilde{c})\\
(a,b,c)  \mapsto  e_2\\
\end{array}
\begin{array} {l}
 \mathbf h_4(v) = id(v)  \\
0  \mapsto  0\\
e_1  \mapsto  e_1\\
e_2  \mapsto  e_2\\
(a,b,c)  \mapsto  (a,b,c)\\
\end{array}
\end{equation}

For the hypothesised maps, we must have $\tilde{c} = \pm c$, for 
unimodular maps preserve volume.  We note that $vol[T(a,b,c)] = |c|/6$ and $vol(T(\tilde{a}, \tilde{b}, \tilde{c})) = |\tilde{c}|/6$,
therefore, $|c| = |\tilde{c}|$.  Without loss of generality, $c = \tilde{c}$ since $(x,y,z) \mapsto (x, y, -z)$
is unimodular.

Also, the maps force certain congruences on $\tilde{a}$  and $\tilde{b}$, specifically:

\begin{equation}
\begin{array}{c}
\mathbf h_1(v)\\
\tilde{a} \equiv (a+b-1)^{-1} a \bmod c\\
\tilde{b} \equiv (a+b-1)^{-1} b \bmod c\\
\end{array}
\begin{array}{c}
\bold h_2(v)\\
\tilde{a} \equiv a^{-1} \bmod c\\
\tilde{b} \equiv -a^{-1} b \bmod c\\
\end{array}
\begin{array}{c}
\bold h_3(v)\\
\tilde{a} \equiv -b^{-1} a \bmod c\\
\tilde{b} \equiv b^{-1} \bmod c\\
\end{array}
\end{equation}

\noindent
How these congruences were derived will be explained via a sample calculation for $\bold h_3(v)$, as they are all similar.  Label
$\bold h_3(v) = A v + b$, and we note that $\bold h_3(0) = 0$.  Then $A e_1 = e_1$ determines the first column of $A$.  
Similarly, $A e_2 = \tilde{a} e_1 + \tilde{b} e_2 + \tilde{c} e_3$, provides the second column.  Column three is determined by 
$A e_3$.  To compute this, we look at the action of $A$ on $a e_1 + b e_2 + c e_3$. 
\begin{equation}
\begin{array}{lcl}
A (a e_1 + b e_2 +c e_3) &  = & a A e_1 + b A e_2 + c A e_3  \\
 &  = & a e_1 + b (\tilde{a} e_1 + \tilde{b} e_2 + \tilde{c} e_3) + c A e_3  \\
 & = & (a + b \tilde{a})e_1 + b \tilde{b} e_2 +b \tilde{c} e_3 + c A e_3 \\
 & = & e_2
\end{array}
\end{equation}

\noindent
Recalling from the above arguement that $c = \tilde{c}$, we find $A e_3$ is:

\begin{equation}
A e_3 = - \frac{(a + b \tilde{a})}{c} e_1 + \frac{(1-b \tilde{b})}{c} e_2 - b e_3
\end{equation}

\noindent
We have now found explicitly the matirx $A$ with $det(A) = 1$.
\begin{equation}
A =
\left[
\begin{array}{ccc}
| & | & | \\
A e_1 & A e_2 & A e_3 \\
| & | & | \\
\end{array}
\right] 
=
\left[
\begin{array}{ccc}
1 & \tilde{a} &  - \frac{(a + b \tilde{a})}{c}\\
0 & \tilde{b} &  \frac{(1-b \tilde{b})}{c}\\ 
0 & c & -b\\
\end{array}
\right]
\end{equation}

\noindent
We recall that $A = A_{ij}$ is an element of $GL(3,{\Z})$, and therefore require that $A_{13} \; \mbox{and} \;  A_{23}$ 
be integers.  This, therefore, forces congruence relations on $\tilde{a}  \; \mbox{and} \; \tilde{b}$.  The congruence relations are as follows:
\begin{equation}
\begin{array}{c}
\tilde{a} \equiv - b^{-1}a \bmod c \ \ \ \ \tilde{b} \equiv b^{-1} \bmod c
\end{array}
\end{equation}

Finally, note all congruences of (\ref{Eq:maps}) may be solved by Fact 1 of the introduction.  Thus, maps $\bold h_1(v)$, $\bold h_2(v)$,
$\bold h_3(v)$, and $\bold h_4(v)$ exist.

\subsection{Naming Cases For Sub-Tetrahedra}

We now define the notation that will allow us to search for such lattice barycentric tetrahedra.  We begin by defining the sub-tetrahedra,
$\mathcal{T}_j$.  We do this so that the $j^{th}$ sub-tetrahedron is brought into the ground position by map $\bold h_j(v)$.
\begin{equation}
\begin{array}{ll}
\mathcal{T}_1 = T(e_1, e_2, BC[T(a,b,c)], (a,b,c))& \\
\mathcal{T}_2 = T(0, e_2, BC[T(a,b,c)], (a,b,c))& \\
\mathcal{T}_3 = T(0, e_1, BC[T(a,b,c)], (a,b,c))& \\
\mathcal{T}_4 = T(\alpha, \beta, \gamma)& \\
\end{array}
\end{equation}

We now define the notation for primitivity.  Any tetrahedron $T(l,m,n)$ can, from Fact 2, be primitive in three different ways.
They are $l \equiv 1 \bmod \gamma$, or $m \equiv 1 \bmod \gamma$, or $l+m \equiv 0 \bmod \gamma$.  These conditions are labeled cases a, b, c respectively.
Thus saying a tetrahedron is lattice barycentric by case (1a, 2c, 3b, 4a) means that \[e_1 \cdot \bold h_1(\alpha, \beta, \gamma) \equiv 1 \bmod \gamma,\;
e_1 \cdot \bold h_2(\alpha, \beta, \gamma) + e_2 \cdot \bold h_2(\alpha, \beta, \gamma) \equiv 0 \bmod \gamma\] \[e_2 \cdot \bold h_3(\alpha, \beta, 
\gamma) \equiv 1 \bmod \gamma,\; \alpha \equiv 1 \bmod \gamma.\]

\section{Analysis in the case $a \equiv 3 \bmod c$ or $b \equiv 3 \bmod c$}\label{Sec:3a}
When we begin to look for such possible lattice barycentric tetrahedra, $T(a,b,c)$, the following lemma greatly reduces the set of possible configurations.
The author thanks his R.E.U. adviser, Stephen Bullock, for this observation.

\begin{lemma} 
\label{lem:8NoDivc}
If $T(a,b,c)$ is lattice barycentric, then $4 \vert c$, however, $8 \nmid c$.
\end{lemma}

\begin{proof} Assume by way of contradiction, that $T(a,b,8\tau), \tau \in {\Z}$ is lattice barycentric.  
Then $T(\frac {a+1} {4}, \frac {b+1} {4}, 2\tau)$ is primitive, and in particular fundamental.   Therefore by Fact 1, 
$a \equiv 3 \bmod 8$, else $GCD(\frac {a+1}{4}, 2\gamma)$ is even, implying boundry points on $T(\frac{a+1}{4}, \frac{b+1}{4}, 2\tau)$.  
Similarly for $b \equiv 3 \bmod 8$.

Now, transform unimodularly $a \rightarrow -a^{-1}b \bmod 8$ and $b \rightarrow b^{-1} \bmod 8$ by the map $\bold{h}_2(v)$ of equation \ref{Eq:maps} .  
Note that $a^{-1} = 3 \bmod 8$, so $\tilde{a} = -3 \cdot 3 = - 9 = 7 \bmod 8$.  Contradiction. \end{proof}

With the above lemma, we begin to look at all such cases of the form $a \equiv 3 \bmod c$.  We note that looking at $b \equiv 3 \bmod c$ is equivalent, as
$(x,y,z) \mapsto (y,x,z)$ is unimodular.  We thus seek configurations of the form $T(3, b, 4 \gamma)$.  Our computation will show
that $a \equiv 3 \bmod c$ or $b \equiv 3 \bmod c$ implies $\gamma = 1$, so $T$ is $T(3,3,4)$.  Otherwise, we would arrive at inconsistencies in the 
congruences.  This section allows us to ignore the case of $a \equiv 3 \bmod c$ or $b \equiv 3 \bmod c$, in the more general section \ref{Sec:gen}, that
is, 4a and 4b.

For $T(3,b,4 \gamma)$, each case of (*,*,*,4a) implies a certain congruence on b.  We provide an example of the congruence which results on $b$ from the 
case 2a of Table \ref{Tab:cona3}.  Label $\bold h_2(\alpha , \beta, \gamma) = (\tilde{\alpha}, \tilde{\beta}, \gamma)$.

\begin{table}

\begin{center}
\begin{tabular}{|c|c|}
\hline
Case 1a & $b \equiv -1 \bmod \gamma$ \\
\hline
Case 1b & $b \equiv -3 \bmod \gamma$ \\
\hline
Case 1c & $3b \equiv -7 \bmod \gamma$ \\
\hline
Case 2a & $1 \equiv 9 \bmod \gamma \Rightarrow \gamma = 8$ \\
\hline
Case 2b & $b \equiv -9 \bmod \gamma$ \\
\hline
Case 2c & $b \equiv 7 \bmod \gamma$ \\
\hline
Case 3a & $b \equiv -1 \bmod \gamma$ \\
\hline
Case 3b & $3b \equiv 1 \bmod \gamma$ \\
\hline
Case 3c & $b \equiv 1 \bmod \gamma$ \\
\hline
\end{tabular}
\end{center}

\caption{\label{Tab:cona3}Congruences with $a \equiv 3 \bmod c$}
\end{table}

\begin{equation}
\begin{array}{rcl}
\tilde{\alpha} \equiv 1 \bmod \gamma & \Longrightarrow & \frac{3^{-1} + 1}{4} \equiv 1 \bmod \gamma  \\
 & \Longrightarrow & 3^{-1} = 1 \equiv 4 \bmod \gamma \\
 & \Longrightarrow & 3^{-1} \equiv 3 \bmod \gamma \\
 & \Longrightarrow & 1 \equiv 9 \bmod \gamma\\
\end{array}
\end{equation}

As an aside, the congruence for Case 2a is a contradiction by lemma \ref{lem:8NoDivc}, given 4a..  We shall
denote these cases by a $\star$.  Whilst it is implied that $\gamma = 8$ by case 2a, we shall say explicitly for the other cases that $\gamma = 8$,  
as they result from two non-obvious congruences.   The results of this exercise are found in Table \ref{Tab:a3res}.  The conclusion of the case study is 
that any lattice barycentric tetrahedron, $T(3,b,4 \gamma)$ is in fact $T(3,3,4)$.

\begin{table}

\begin{center}
\begin{tabular}{|c|c|}
\hline
Configuration & Conclusion \\
\hline
1a, 2a, 3a, 4a & $\star$ \\
\hline
1a, 2a, 3b, 4a & $\star$ \\
\hline
1a, 2a, 3c, 4a & $\star$ \\
\hline
1a, 2b, 3a, 4a & $\gamma = 8$ \\
\hline
1a, 2b, 3b, 4a & $\gamma = 8$ \\
\hline
1a, 2b, 3c, 4a & $\gamma = 8$ \\
\hline
1a, 2c, 3a, 4a & $\gamma = 8$ \\
\hline
1a, 2c, 3b, 4a & $\gamma = 8$ \\
\hline
1a, 2c, 3c, 4a & $\gamma = 8$ \\
\hline
1b, 2a, 3a, 4a & $\star$ \\
\hline
1b, 2a, 3b, 4a & $\star$ \\
\hline
1b, 2a, 3c, 4a & $\star$ \\
\hline
1b, 2b, 3a, 4a & $\gamma = 8$ \\
\hline
1b, 2b, 3b, 4a & $\gamma = 3 \longrightarrow  b \equiv 0 \bmod 3$ and $b \equiv 3^{-1} \bmod 3$  \\
\hline
1b, 2b, 3c, 4a & $\gamma = 8$ \\
\hline
1b, 2c, 3a, 4a & $\gamma = 8$ \\
\hline
1b, 2c, 3b, 4a & $\gamma = 5$ \\
\hline
1b, 2c, 3c, 4a & $\gamma = 3$ or $\gamma = 5$ \\
\hline
1c, 2a, 3a, 4a & $\star$ \\
\hline
1c, 2a, 3b, 4a & $\star$ \\
\hline
1c, 2a, 3c, 4a & $\star$ \\
\hline
1c, 2b, 3a, 4a & $\gamma = 8$ \\
\hline
1c, 2b, 3b, 4a & $\gamma = 8$ \\
\hline
1c, 2b, 3c, 4a & $\gamma = 5$ \\
\hline
1c, 2c, 3a, 4a & $\gamma = 8$ \\
\hline
1c, 2c, 3b, 4a & $\gamma = 8$ \\
\hline
1c, 2c, 3c, 4a & $\gamma = 6$ \\
\hline
\end{tabular}
\end{center}

\caption{\label{Tab:a3res}Case study for $a \equiv 3 \bmod c$}
\end{table}

\section{Remaining Cases}\label{Sec:gen}
\subsection{Congruences For Primitive Sub-Tetrahedra}

We shall now look at a sample calculation for the conditions on primitivity and the resulting congruences on $a$ and $b$ in the general case
$a \not\equiv 3 \bmod c$ and $b \not\equiv 3 \bmod c$.  Let us take for example, 
the map $\bold h_1(v)$.  Map $\bold h_1(v)$ takes $BC[T(a,b,c)] = (\alpha, \beta, \gamma) \mapsto (\tilde{\alpha}, \tilde{\beta}, \tilde{\gamma})$.  
We know that $\tilde{\alpha} \equiv \frac{\tilde{a}+1}{4} \bmod \gamma$ and $\tilde{\beta} \equiv 
\frac{\tilde{b}+1}{4} \bmod \gamma$.  We shall look at the first condition for primitiveness 1a, as 1b and 1c are similar.

\begin{equation}
\begin{array}{lcl}
\tilde{\alpha} \equiv 1 \bmod \gamma & \Longrightarrow & \frac{(a+b-1)^{-1} a + 1}{4} \equiv 1 \bmod \gamma  \\
 & \Longrightarrow & (a+b-1)^{-1} a + 1 \equiv 4 \bmod \gamma \\
 & \Longrightarrow & a \equiv 3 (a+b-1) \bmod \gamma \\
 & \Longrightarrow & 2a + 3b \equiv 3 \bmod \gamma\\
\end{array}
\end{equation}

\noindent
Upon checking all cases, we arrive at Table \ref{Tab:con}.
 
\begin{table}
   
\begin{center}
\begin{tabular}{|c|c|}
\hline
Case 1a & $2a + 3b \equiv 3 \bmod \gamma$ \\
\hline
Case 1b & $3a + 2b \equiv 3 \bmod \gamma$ \\
\hline
Case 1c & $3a + 3b \equiv 2 \bmod \gamma$ \\
\hline
Case 2a & $3a  \equiv 1 \bmod \gamma$ \\
\hline
Case 2b & $3a + b \equiv 0 \bmod \gamma$ \\
\hline
Case 2c & $2a - b \equiv -1 \bmod \gamma$ \\
\hline
Case 3a & $a + 3b \equiv 0 \bmod \gamma$ \\
\hline
Case 3b & $3b \equiv 1 \bmod \gamma$ \\
\hline
Case 3c & $a - 2b \equiv 1 \bmod \gamma$ \\
\hline
Case 4a & $a \equiv 3 \bmod \gamma$ \\
\hline
Case 4b & $b \equiv \bmod \gamma$ \\
\hline 
Case 4c & $ a + b \equiv -2 \bmod \gamma$ \\
\hline
\end{tabular}
\end{center} 

\caption{\label{Tab:con}General congruences on $a$ and $b$}
\end{table}

\subsection{Remaining Case Study With Examples}
We will now look at the general case, $a \not\equiv 3 \bmod c$ and $b \not\equiv 3 \bmod c$.  Recall section \ref{Sec:3a} showed $T(3,3,4)$ is the only
lattice barycentric tetrahedron $T$, of the form $T(3, b, 4 \gamma)$.  Two sample calculations from the generic case study will now be shown.  
The first will use a configuration in which a 
unimodular equivalence class is found, the other will use a configuration that leads to an inconsistency.

\subsubsection*{Example Of Calculation Resulting In A Unimodular Equivalence Class}

\noindent
Let us look at the case 1b, 2a, 3a, 4c, from Table \ref{Tab:con}.  This forces the following congruences:

\[
\begin{array}{cccc}
3a + 2b \equiv 3 \bmod \gamma & 3a \equiv 1 \bmod \gamma & a + 3b \equiv 0 \bmod \gamma & a +b \equiv -2 \bmod \gamma
\end{array}
\]

One may solve this system of linear congruences by row reductions.
However, since we seek congruence relations on $a$ and $b$ modulo
$\gamma = \frac{c} {4}$, we do not divide by any integers except $2$.
Division by $2$ is allowed modulo $\gamma$ by lemma \ref{lem:8NoDivc}.  This point in the reduction is indicated by {\bf !!}.
\begin{displaymath}
\left[
\begin{array}{cc|c}
3 & 2 & 3 \\
3 & 0 & 1 \\
1 & 3 & 0 \\
1 & 1 & -2 \\ 
\end{array}
\right]
\longrightarrow
\left[
\begin{array}{cc|c}
0 & 0 & 0 \\
0 & -3 & 7 \\
0 & 2 & 2 \\
1 & 1 & -2 \\ 
\end{array}
\right]
{\buildrel {\bf !!} \over \longrightarrow}
\left[
\begin{array}{cc|c}
0 & 0 & 0 \\
0 & -3 & 7 \\
0 & 1 & 1 \\
1 & 1 & -2 \\ 
\end{array}
\right]
\longrightarrow
\left[
\begin{array}{cc|c}
0 & 0 & 0 \\ 
0 & 0 & 10 \\
0 & 1 & 1 \\
1 & 0 & -3 \\ 
\end{array}
\right]
\end{displaymath}

We see that the second line demands that $10 \equiv 0 \bmod \gamma$, so $\gamma = 0,2,5,10$.  By lemma \ref{lem:8NoDivc} and $\gamma > 0$, 
$\gamma = 5$.  Then by Chinese Remainder Theorem, $a \equiv -3 \bmod \gamma$ and $a \equiv 3 \bmod 4$ means $a \equiv 7 \bmod 20$.  Similarly, 
by CRT $b \equiv 1 \bmod 5$ and $b \equiv 3 \bmod 4$ means $b \equiv 11 \bmod 20$.  We note also that $(7, 11, 20)$ satisfies the conditions for
primitivity and fundamentality, so this case produces $T(7,11,20)$.

\subsubsection*{Example Of Calculation Not Resulting In A Unimodular Equivalence Class}

Let us look now at the case 1a, 2a, 3a, 4c from Table \ref{Tab:con}.  This forces the following congruences:

\[
\begin{array}{cccc}
2a + 3b \equiv 3 \bmod \gamma & 3a \equiv 1 \bmod \gamma & a + 3b \equiv 0 \bmod \gamma & a + b \equiv -2 \bmod \gamma
\end{array}
\]

We apply the same solution technique as in the last example.
\begin{displaymath}
\left[
\begin{array}{cc|c}
2 & 3 & 3 \\
3 & 0 & 1 \\
1 & 3 & 0 \\
1 & 1 & -2 \\ 
\end{array}
\right]
\longrightarrow
\left[
\begin{array}{cc|c}
0 & 1 & 7 \\
0 & -3 & 7 \\
0 & 2 & 2 \\
1 & 1 & -2 \\ 
\end{array}
\right]
{\buildrel {\bf !!} \over \longrightarrow}
\left[
\begin{array}{cc|c}
0 & 1 & 7 \\
0 & 0 & 28 \\
0 & 1 & 1 \\
1 & 1 & -2 \\ 
\end{array}
\right]
\longrightarrow
\left[
\begin{array}{cc|c}
0 & 0 & 2 \\ 
0 & 0 & 0 \\
0 & 1 & 1 \\
1 & 0 & -3 \\ 
\end{array}
\right]
\end{displaymath}

We see that the second line demands that $2 \equiv 0 \bmod \gamma$.  By lemma \ref{lem:8NoDivc} $\gamma$ is odd, and therefore we have an
incomsistency.

\subsection*{Results Of The Case Study}
We now present the results of an exhaustive search for the equivalence classes in Table \ref{Tab:Res}.  The final conclusion is that $\gamma = 1$ or 
$\gamma = 5$, producing $T(3,3,4)$ and $T(7,11,20)$ consecutively.  This concludes the proof of proposition \ref{Prop:main}.

\begin{table}
\begin{center}
\begin{tabular}{|c|c|}
\hline
Combination & Conclusion \\
\hline
1a, 2a, 3a, 4c & $\gamma$ even \\
\hline
1a, 2a, 3b, 4c &  $\gamma$ even \\
\hline
1a, 2a, 3c, 4c & $\gamma$ even \\
\hline
1a, 2b, 3a, 4c & $\gamma$ even \\
\hline
1a, 2b, 3b, 4c & $\gamma = 5 \rightarrow (7, 11, 20)$  \\
\hline
1a, 2b, 3c, 4c & $\gamma$ even \\
\hline
1a, 2c, 3a, 4c & $\gamma = 3 \rightarrow (3, 7, 12) (3,12) \neq 1$ \\
\hline
1a, 2c, 3b, 4c & $\gamma$ even \\
\hline
1a, 2c, 3c, 4c & $\gamma = 3 \rightarrow (3, 7, 12) (3,12) \neq 1$ \\
\hline
1b, 2a, 3a, 4c & $\gamma = 5 \rightarrow (7, 11, 20)$ \\
\hline
1b, 2a, 3b, 4c & $\gamma$ even \\
\hline
1b, 2a, 3c, 4c & $\gamma$ even \\
\hline
1b, 2b, 3a, 4c & $\gamma$ even \\
\hline
1b, 2b, 3b, 4c & $\gamma$ even \\
\hline
1b, 2b, 3c, 4c & $\gamma = 3 \rightarrow (3, 7, 12) (3,12) \neq 1$ \\
\hline
1b, 2c, 3a, 4c & $\gamma$ even \\
\hline
1b, 2c, 3b, 4c & $\gamma$ even \\
\hline
1b, 2c, 3c, 4c & $\gamma = 3 \rightarrow (3, 7, 12) (3,12) \neq 1$ \\
\hline
1c, 2a, 3a, 4c & $\gamma$ even \\
\hline
1c, 2a, 3b, 4c & $\gamma$ even \\
\hline
1c, 2a, 3c, 4c & $\gamma$ even \\
\hline
1c, 2b, 3a, 4c & $\gamma$ even \\
\hline
1c, 2b, 3b, 4c & $\gamma$ even \\
\hline
1c, 2b, 3c, 4c & $\gamma$ even \\
\hline
1c, 2c, 3a, 4c & $\gamma$ even \\
\hline
1c, 2c, 3b, 4c & $\gamma$ even \\
\hline
1c, 2c, 3c, 4c & $\gamma$ even \\
\hline
\end{tabular}
\end{center} 
\caption{\label{Tab:Res} Results of generic case study}
\end{table}

\section{Acknowledgments}

Many thanks to Dr. Stephen Bullock, formerly of the mathematics department at the University of Michigan, currently employed by M.C.S.D. of N.I.S.T. at
the Gaithersburg, Maryland campus.  I learned a great deal and enjoyed very much our many conversations.  I also thank Dr. Bruce Reznick of the mathematics
department at the University of Illinois - Urbana/Champaign and Dr. Paul Federbush of the mathematics department at the University of Michigan.

\end{document}